\documentclass[11pt,reqno]{amsart}
\usepackage{indentfirst, amssymb,amsmath,amsthm}
\usepackage{enumitem}
\usepackage{xcolor}
\usepackage[colorlinks]{hyperref}
\usepackage{tikz-cd}
\newtheorem{theorem1}{Theorem}[section]
\newtheorem{example1}[theorem1]{Example}
\newtheorem{lemma1}[theorem1]{Lemma}

\newtheorem{corollary1}[theorem1]{Corollary}
\newtheorem{remark1}[theorem1]{Remark}

\newtheorem{definition1}[theorem1]{Definition}

\newtheorem{note1}[theorem1]{Note}

\newcommand{\be}{\begin{equation}}
\newcommand{\ee}{\end{equation}}
\newcommand{\beas}{\begin{eqnarray*}}
\newcommand{\eeas}{\end{eqnarray*}}
\newcommand{\bea}{\begin{eqnarray}}
\newcommand{\eea}{\end{eqnarray}}

\numberwithin{equation}{section}

\begin{document}

\setlength{\unitlength}{1mm} \baselineskip .52cm
\setcounter{page}{1}
\pagenumbering{arabic}
\title[On a generalization of quasi-metric space]{On a generalization of quasi-metric space}

\author[Sugata Adhya and A. Deb Ray]{Sugata Adhya and A. Deb Ray}

\address{Department of Mathematics, The Bhawanipur Education Society College. 5, Lala Lajpat Rai Sarani, Kolkata 700020, West Bengal, India.}
\email {sugataadhya@yahoo.com}

\address{Department of Pure Mathematics, University of Calcutta. 35, Ballygunge Circular Road, Kolkata 700019, West Bengal, India.}
\email {debrayatasi@gmail.com}

\maketitle

\begin{abstract}
We find an extension of the quasi-metric (to be called
$g$-quasi metric) such that the induced generalized topology may fail to form a topology. We show that $g$-quasi metrizability is a $g$-topologically invariant property of generalized topological spaces. Extending metric product and uniform continuity for $g$-quasi metric spaces, we note that a $g$-quasi metric may fail to be uniformly continuous in the extended sense unlike usual metric. Finally, we extend the study of completeness, Lebesgue property and weak $G$-completeness for $g$-quasi metric spaces.

\end{abstract}


\section{\textbf{Introduction}}

Cs{\'a}sz{\'a}r \cite{czaszar1} proposed a notion of generalized topology by taking into account the idea of monotone mappings. It accommodates various open-like sets existed in literature \cite{csaszar1997generalized,levine1963semi,abd1982precontinuous,Njstad1965OnSC}. Given a nonempty set $X,$ it is defined as  a subcollection of $\mathcal P(X)$ which contains $\emptyset$ and is closed under arbitrary union. Considering the members of a generalized topology as open, it then became natural to study the usual topological notions for generalized topology. Accordingly, analogoues of closed set, closure, interior, product and subspaces, continuous functions, countability and separation axioms, compactness and connectedness have been studied for generalized topological spaces. Additionally, multiple weaker versions of the above notions have also been investigated in the light of weaker forms of open sets in the context of generalized topology. The interested readers may consult \cite{wu-zhu} and references therein.

In view of the facts that the associated open balls in any metric space forms a base for a natural topology and generalized topology is a generalization of usual topology, it is natural to search for a generalization of metric structure that under standard approach induces a generalized topology and fails to produce a topology, in general. This paper addressed this question. 

Here, in Section \ref{section 3}, we obtain the related spaces (termed as $g$-quasi metric spaces) by generalizing Wilson’s widely studied notion of the quasi-metric structure \cite{wilson1931quasi}. Subsequently, we discuss certain separation properties of the generalized topology induced by a $g$-quasi metric. We demonstrate that $g$-quasi metrizability is an invariant property of the generalized topology. Next, in Section \ref{section 4}, we propose the natural extensions of the notions of metric product and uniform continuity in the context of $g$-quasi metric spaces. It is noted that, unlike usual metric, a $g$-quasi metric may fail to be uniformly continuous in the extended sense while considered as a mapping from the product space to $\mathbb R$. 

The study made in Section \ref{section 3} and Section \ref{section 4} pave the way for the related investigations on Cauchyness and completeness. Accordingly, in Section \ref{section 5}, we take up the task of extending them in $g$-quasi metric spaces. Apart from the usual completeness, we introduce two stronger forms of completenesses viz. Lebesgue property \cite{l1,l2,l3} and weak $G$-completeness \cite{1sa-g-com,2greg-g-com,3greg2g-com} in $g$-quasi metric spaces. It is known that both Lebesgue property and weak $G$-completeness are intermediate metric properties between compactness and completeness that can be characterized in terms of pseudo-Cauchy \cite{l3} and $G$-Cauchy \cite{2greg-g-com} sequences respectively. In what follows, we explore the mutual dependence of those completenesses for $g$-quasi metrics and enquire their behaviours in the product spaces through Cauchy, pseudo-Cauchy and $G$-Cauchy sequences.

\section{\textbf{Preliminaries}}

This section discusses the prerequisites that will be required subsequently. 

\begin{definition1}
\normalfont\cite{czaszar1,csaszar2004extremally,mashhour1983supratopological} A generalized topology\index{generalized!topology} $\mu$ on a nonempty set $X$ is a collection of its subsets such that $\emptyset\in\mu$ and $\mu$ is closed under arbitrary union. The pair $(X,\mu)$ is called a generalized topological space. Moreover, if $X\in\mu$ then $\mu$ is called a supratopology\index{supratopology} or a strong generalized topology\index{strong generalized topology}.

Given a generalized topological space $(X,\mu),$ the elements of $\mu$ are called generalized open sets\index{generalized!open set} (or $\mu$-open sets\index{$\mu$-open set}) and their complements are called generalized closed sets\index{generalized!closed set} (or $\mu$-closed sets\index{$\mu$-closed set}) in $(X,\mu)$.   
\end{definition1}

\begin{definition1}
\normalfont\cite{csaszar2007modification} Given a nonempty set $X$ and $\mathcal{B}\subset\mathcal{P}(X)$ with $\emptyset\in\mathcal{B},$ all possible unions of elements of $\mathcal{B}$ form a generalized topology $\mu(\mathcal{B})$ on $X.$ Here $\mathcal{B}$ is called a base\index{base for!a generalized topology} for $\mu(\mathcal{B}).$ Equivalently, $\mu(\mathcal{B})$ is said to be generated by the base $\mathcal{B}.$
\end{definition1}

\begin{definition1}
\normalfont\cite{czaszar1,czaszarproduct} Let $(X,\mu)$ and $(Y,\mu')$ be two generalized topological spaces. 

(a) A mapping $f: X\to Y$ is called generalized continuous\index{generalized!continuous mapping} or $(\mu,\mu')$-continuous\index{$(\mu,\mu')$-continuous mapping} if $f^{-1}(G)\in\mu,~\forall~G\in\mu'$. 

(b) A mapping $f:X\to Y$ is called a generalized homeomorphism\index{generalized!homeomorphism} or $(\mu,\mu')$-homeomorphism\index{$(\mu,\mu')$-homeomorphism} if $f$ is bijective and $f,f^{-1}$ are generalized continuous. 

(c) A property of generalized topological spaces that is invariant under generalized homeomorphism is called a $g$-topologically invariance\index{g@$g$-topologically invariance}.
\end{definition1}

\begin{definition1}
\normalfont\cite{sarsak_weak_sep} Let $(X,\mu)$ be a strong generalized topological space.

(a) $(X,\mu)$ is called $\mu$-$T_0$\index{$\mu$-$T_0$} if for $x,y\in X$ with $x\ne y$ there exists $B\in\mu$ such that exactly one of $x$ and $y$ is in $B.$

(b) $(X,\mu)$ is called $\mu$-$T_1$\index{$\mu$-$T_1$} if for $x,y\in X$ with $x\ne y$ there exist $B_1,B_2\in\mu$ such that $B_1$ contains $x$ but not $y$ and $B_2$ contains $y$ but not $x.$

Clearly every $\mu$-$T_1$ strong generalized topological space is $\mu$-$T_0.$
\end{definition1}

\begin{theorem1}
\normalfont\cite{sarsak_weak_sep} Every singleton set in a $\mu$-$T_1$ strong generalized topological space is $\mu$-closed.
\end{theorem1}

\begin{definition1}
\normalfont\cite{raybhowmik2015} Let $X$ be a nonempty set. A nonempty subset $\mathcal U$ of $\mathcal{P}(X\times X)$ is called a generalized quasi-uniformity (or $g$-quasi uniformity) if 

(i) each member of $\mathcal U$ contains the diagonal of $X,$ 

(ii) $\mathcal U$ is closed under superset, 

(iii) for $U\in\mathcal U$ there exists $V\in\mathcal U$ such that $V\circ V\subset U.$

\noindent In this case, $(X,\mathcal U)$ is called a generalized quasi-uniform space (or $g$-quasi uniform space).
\end{definition1}

\begin{theorem1}
\normalfont\cite{raybhowmik2015} Given a nonempty set $X$ and $\mathcal{B}~(\ne\emptyset)\subset\mathcal P(X\times X),$ $\mathcal{B}$ forms a base for some $g$-quasi uniformity on $X$ if and only if (i) $\Delta(X)\subset B,~\forall~B\in\mathcal B,$ and (ii) $B\in\mathcal B\implies\exists~V\in\mathcal B$ such that $V\circ V\subset B.$  

Moreover, such $\mathcal{B}$ is a base for the $g$-quasi uniformity $\{V\subset X:B\subset V\text{ for some }B\in\mathcal{B}\}$ on $X.$
\end{theorem1}

we finish this section by recalling certain preliminaries on Lebesgue property and weak $G$-completeness for metric spaces.

\begin{definition1}
\normalfont\cite{l2} A metric space on which every real-valued continuos function is uniformly continuous is said to be Lebesgue (or Atsuji space). 
\end{definition1}

\begin{definition1}
\normalfont\cite{l2} A sequence $(x_n)$ in a metric space $(X,d)$ is said to be pseudo-Cauchy if given $\epsilon>0,k\in\mathbb N$ there exists distinct $m,n~(>k)\in\mathbb N$ such that $d(x_m,x_n)<\epsilon.$
\end{definition1}

\begin{theorem1}
\normalfont\cite{l2,l3} A metric space is Lebesgue if and only if every pseudo-Cauchy sequence having distinct terms clusters in it. 
\end{theorem1}

\begin{definition1}
\normalfont\cite{1sa-g-com,2greg-g-com,3greg2g-com} A sequence $(x_n)$ in a metric space $(X,d)$ is called $G$-Cauchy\index{g@$G$-Cauchy sequence!in metric spaces} if $\lim\limits_{n\to\infty}d(x_{n+p},x_n)=0,~\forall~ p\in\mathbb N$ (or equivalently, $\lim\limits_{n\to\infty}d(x_{n+1},x_n)=0$). A metric space in which every $G$-Cauchy sequence converges is said to be weak $G$-complete.
\end{definition1}

Both Lebesgue property and weak $G$-completeness are strictly intermediate between compactness and completeness of metric spaces \cite{2greg-g-com,l2}.

\section{\textbf{\texorpdfstring{$g$}--Quasi Metric Spaces and the Induced Generalized Topology}}\label{section 3}

\begin{definition1}
\normalfont\cite{wilson1931quasi} Let $X$ be a nonempty set. A mapping $d:X\times X\to\mathbb R$ is called a quasi-metric on $X,$ if:

(a) $\forall~x,y\in X,~d(x,y)\ge 0$ and $d(x,y)=0\iff x=y;$

(b) $\forall~x,y,z\in X,~d(x,y)\le d(x,z)+d(z,y).$

Here $d$ is called a quasi-metric\index{quasi-metric (space)} on $X.$ Moreover, the pair $(X,d)$ is called a quasi-metric space.
\end{definition1}

\begin{definition1}
\normalfont Let $X$ be a nonempty set. A mapping $d:X\times X\to\mathbb R$ is called a $g$-quasi metric\index{g@$g$-quasi!metric (space)} on $X,$ if there exists $r\ge0$ such that

(a) $\forall~x,y\in X,~d(x,y)\ge r$ and $d(x,y)=r\iff x=y;$

(b) $\forall~x,y,z\in X,~d(x,y)\le d(x,z)+d(z,y).$

Here $d$ is called a $g$-quasi metric on $X$ and $r,$ the index\index{index of $g$-quasi metric (space)} of $d.$ Moreover, the pair $(X,d)$ is called a $g$-quasi metric space (of index $r).$

Clearly a $g$-quasi metric of index $0$ is a quasi-metric and vice versa.
\end{definition1}

\begin{definition1}\normalfont
A $g$-quasi metric $d$ on a nonempty set $X$ (and hence the related $g$-quasi metric space $(X,d)$) is said to be symmetric\index{symmetric!g@$g$-quasi metric (space)} if $d(x,y)=d(y,x),~\forall~x,y\in X.$

Clearly a symmetric $g$-quasi metric of index $0$ is a metric and vice versa.
\end{definition1}

\begin{theorem1}
\normalfont Let $(X,d)$ be a quasi-metric space. Then for $r\ge0,~d'=d+r$ forms a $g$-quasi metric on $X$ of index $r.$
\end{theorem1}

\begin{proof}
\normalfont (a) Clearly $\forall~x,y\in X,~d'(x,y)\ge r$ and $d'(x,y)=r\iff x=y;$

(b) Choose $x,y,z\in X.$ Then $d'(x,y)=d(x,y)+r\le d(x,z)+d(z,y)+r\le\{d(x,z)+r\}+\{d(z,y)+r\}\le d'(x,z)+d'(z,y).$ 

Hence the result follows.
\end{proof}

However given a $g$-quasi metric $d$ on a nonempty set $X,$ $d'=d-r$ may not form a quasi-metric on it for all choices of $r\ge0:$

\begin{example1}
\normalfont Let $X=[2,4]$ and $d:X\times X\to\mathbb R$ be given by $d(x,y)=(x-y)^2+100,~\forall~x,y\in X.$ Then for no values of $r\ge0,$ $d'=d-r$ forms a quasi-metric on $X.$ It follows by observing that for $r=100,~d'(2,4)>d'(2,3)+d'(3,4),$ while for all other values of $r,~d'(x,x)\ne0,~\forall~x\in X.$ 

However $d$ is a $g$-quasi metric on $X$ of index $100:$

(a) $\forall~x,y\in X,~d(x,y)\ge 100$ and $d(x,y)=100\iff x=y,$

(b) $\forall~x,y,z\in X,~d(x,y)+d(y,z)\ge100+100\ge(x-z)^2+100=d(x,z).$
\end{example1}

\begin{definition1}
\normalfont Let $(X,d)$ be a $g$-quasi metric space of index $r\ge0.$ Given $x\in X$ and $p>0,$ we denote the set $\{y\in X:d(x,y)<p\}$ by $B_d(x,p)$ (or simply by $B(x,p)$). Clearly $\mathcal B(d)=\{B(x,p):x\in X,p>0\}\bigcup\{\emptyset\}$ forms a base for some strong generalized topology $\mu(d)$ on $X.$ It is called the generalized topology\index{generalized!topology induced by a $g$-quasi metric} induced by $d.$
\end{definition1}

\begin{note1}\normalfont
It should be noted, at this stage, that if $(X,d)$ is a symmetric $g$-quasi metric space of index $0$ (i.e., a metric space) then $\mathcal{B}(d),$ defined as before, forms a base for the topology induced by $d.$
\end{note1}

In what follows, we show that for all positive values of $r,$ a symmetric $g$-quasi metric space $(X,d)$ can be found so that (i) $d$ is of index $r,$ (ii) $\mu(d)$ does not form a topology on $X.$ We consider the following example.

\begin{example1}\label{gquasi_mainexm}
\normalfont Let $r>0$ and $d:\mathbb R\times\mathbb R\to\mathbb R$ be defined by $$d(x,y)=\begin{cases}
r  & \text{if $x=y$} \\
2r & \text{if $0<|x-y|\le r$}\\
|x-y| & \text{if $|x-y|> r$}
\end{cases}$$

We show that $(\mathbb R,d)$ is a $g$-quasi metric space of index $r$ though $\mu(d)$ does not form a topology on $\mathbb R.$

(a) Clearly $d(x,y)\ge r$ and $d(x,y)=r\iff x=y,$ $\forall~x,y\in\mathbb R.$

(b) Let $x,y,z\in\mathbb R.$ We show that $d(x,y)\le d(x,z)+d(z,y).$ 

Note for $0\le|x-y|\le r$ the above inequality follows from (a). So let us assume $|x-y|>r.$ 

If $|x-z| = 0$ or $|z-y| = 0$ then the inequality is immediate.

If $|x-z|,|z-y|>r$ then it follows from the order property of $\mathbb R.$

If $0<|x-z|,|z-y|\le r$ then $d(x,y)=|x-y|\le|x-z|+|z-y|\le4r=d(x,z)+d(z,y).$

If $|x-z|>r$ and $0<|z-y|\le r$ then $d(x, y)\le d(x,z)+|z-y|\le d(x,z)+r\le d(x,z)+2r=d(x,z)+d(z,y).$

If $|z-y|>r$ and $0<|x-z|\le r$ then it follows similarly as before.

Thus $d$ forms a $g$-quasi metric on $\mathbb R$ of index $r.$
\vspace{1mm}

We now show that $\mu(d)$ does not form a topology on $\mathbb R.$ 

Since $B\left(r,2r+\frac{r}{10}\right)\bigcap $ $B\left(\frac{21r}{5},2r+\frac{r}{10}\right)$ $=\left(r-\frac{21r}{10},r+\frac{21r}{10}\right)\bigcap\left(\frac{21r}{5}-\frac{21r}{10},\frac{21r}{5}+\frac{21r}{10}\right)$ $=\left(\frac{21r}{10},r+\frac{21r}{10}\right),$ it suffices to show that $(\frac{21r}{10},r+\frac{21r}{10})$ does not contain any nonempty set of the form $B(x,c)$ where $x\in\mathbb R,c>0.$

Suppose otherwise. Then $B(x,c)\subset \left(\frac{21r}{10},r+\frac{21r}{10}\right)$ for some $x\in\mathbb R$ and $c>r.$

\textit{Case I:} $r<c\le2r.$ Then for chosen $y\in\mathbb R$ with $r<|x-y|<c,~y\in B(x,c)$ Also $x\in B(x, c).$ 
Thus $x,y\in B(x, c)\subset \left(\frac{21r}{10},r+\frac{21r}{10}\right),$ a contradiction to $|x-y|>r.$

\textit{Case II:} $c>2r.$ Then for each $y\in (x-r,x+r),~d(x,y)<c.$ Consequently $(x-r,x+r)\subset B(x,c)\implies(x-r,x+r)\subset \left(\frac{21r}{10},r+\frac{21r}{10}\right),$ a contradiction to $\left|\left(\frac{21r}{10},r+\frac{21r}{10}\right)\right|=r.$

The contradictions arrived at both the cases prove our claim.
\end{example1}

\begin{remark1}\normalfont
Let $(X,\overline{d})$ be a $g$-quasi metric space with index $r.$ For $\epsilon>r,$ set $V_\epsilon=\{(x,y)\in X\times X:\overline{d}(x,y)<\epsilon\}.$

If $r=0,$ then it is clear that $\mathcal{B}_{\mathcal{U}}=\{V_\epsilon:\epsilon>r\}$ forms a base for a $g$-quasi uniformity on $X$ alike the classical case.

However for all other non-negative values of $r,$ there is some $g$-quasi metric space with index $r$ such that $\mathcal{B}_{\mathcal{U}}$ fails to form a base for some $g$-quasi uniformity on $X$ as we see now.

Consider the $g$-quasi metric space $(\mathbb{R},d),$ as defined in Example \ref{gquasi_mainexm}, with index $r>0.$ If possible, let $\mathcal{B}_{\mathcal{U}}=\{V_\epsilon:\epsilon>r\}$ forms a base for a $g$-quasi uniformity on $\mathbb R$ where $V_\epsilon=\{(x,y)\in X\times X: d(x,y)<\epsilon\},~\forall~\epsilon>r.$ Then there is $\delta>r$ such that, $V_\delta\circ V_\delta\subset V_{\frac{3r}{2}}$.

Set $x=0,y=r+\frac{\delta-r}{2},z=r+\delta.$ Then $d(x,y)=d(y,z)=r+\frac{\delta-r}{2}<\delta$ and hence $(x,y),(y,z)\in V_\delta$ implies $(x,z)\in V_{\frac{3r}{2}}.$ i.e., $d(x,z)<\frac{3r}{2},$ i.e., $r+\delta<\frac{3r}{2}$ i.e., $\delta<\frac{r}{2},$ a contradiction. Hence $\mathcal{B}_{\mathcal{U}}$ is not a base for a $g$-quasi uniformity on $\mathbb R.$
\end{remark1}

\begin{theorem1}
\normalfont Let $d$ be a $g$-quasi metric on $X.$ Then $(X,\mu(d))$ is $\mu$-$T_1$.
\end{theorem1}

\begin{proof}
\normalfont Let $d$ be of index $r.$ Choose $x,y\in X$ such that $x\ne y.$ 

Clearly $d(x,y),~d(y,x)>r.$ Choose $p\in\mathbb R$ such that $r<p<\min\{d(x,y),d(y,x)\}.$ Then $B(x,p),~B(y,p)$ are generalized open sets in $(X,\mu(d))$ containing $x,~y$ respectively such that $x\notin B(y,p)$ and $y\notin B(x,p).$
\end{proof}

Let $d$ be a $g$-quasi metric on $X.$ Then we may conclude the following, stated as corollaries:

\begin{corollary1}
\normalfont  Each singleton set in $(X,\mu(d))$ is $\mu$-closed.
\end{corollary1}

\begin{corollary1}
\normalfont $(X,\mu(d))$ is $\mu$-$T_0$.
\end{corollary1}

\begin{remark1}
\normalfont $g$-Quasi metrices of different indices may induce the same generalized topology. For example, choosing $X=\{x,y\}$ and $d_1,d_2:X\times X\to\mathbb R$ as given by $d_1(x,x)=d_1(y,y)=3,d_1(x,y)=d_1(y,x)=4$ and $d_2(x,x)=d_2(y,y)=5,d_2(x,y)=d_2(y,x)=6,$
we observe that both $d_1,d_2$ induce discrete topology on $X$ though they have distinct indices.
\end{remark1}

\begin{definition1}
\normalfont A generalized topological space $(X,\mu)$ is called $g$-quasi metrizable\index{g@$g$-quasi!metrizable generalized topological space} if for some $g$-quasi metric $d$ on $X,~\mu(d)=\mu.$
\end{definition1}

\begin{theorem1}
\normalfont Let $(X,\mu),(Y,\mu')$ be two generalized topological spaces and $f:X\to Y$ be a generalized homeomorphism. If $(X,\mu)$ is $g$-quasi metrizable, then so is $(Y,\mu').$
\end{theorem1}

\begin{proof}
\normalfont Let $(X,\mu)$ be $g$-quasi metrizable and $d:X\times X\to\mathbb R,$ a $g$-quasi metric on $X$ of index $r\ge0$ such that $\mu(d)=\mu.$

Define $d':Y\times Y\to\mathbb R$ by $d'(y_1,y_2)=d(f^{-1}(y_1),f^{-1}(y_2)),~\forall~y_1,y_2\in Y.$ Then 

(a) $\forall~y_1,y_2\in Y,~d'(y_1,y_2)\ge r$ and $d'(y_1,y_2)=r\iff f^{-1}(y_1)=f^{-1}(y_2)\iff y_1=y_2;$ 

(b) $\forall~y_1,y_2,y_3\in Y,~d'(y_1,y_2)+d'(y_2,y_3)=d(f^{-1}(y_1),f^{-1}(y_2))+d(f^{-1}(y_2),f^{-1}(y_3))\ge d(f^{-1}(y_1),f^{-1}(y_3))= d'(y_1,y_3).$ Thus $d'$ forms a $g$-quasi metric on $Y$ of index $r.$ 

Choose $V\in\mu'$ and $y\in V.$ Then for some $x\in X,~p>r$ we have $f^{-1}(y)\in B_d(x,p)\subset f^{-1}(V)\implies y\in f(B_d(x,p))\subset V\implies y\in B_{d'}(f(x),p)\subset V.$

Thus $\mathcal B(d')$ forms a base for $\mu'.$ Hence the result follows.
\end{proof}

\section{\textbf{Product of \texorpdfstring{$g$}--Quasi Metrics}}\label{section 4}

\begin{theorem1}
\normalfont Let $(X,d_X)$ and $(Y,d_Y)$ be $g$-quasi metric spaces of the same index $r.$ Define $d_{XY}:(X\times Y)\times(X\times Y)\to\mathbb R$ by $$d_{XY}((x_1,y_1),(x_2,y_2))=\max\{d_X(x_1,x_2),d_Y(y_1,y_2)\},$$ $\forall~x_1,x_2\in X$ and $y_1,y_2\in Y.$ Then $d_{XY}$ defines a $g$-quasi metric on $X\times Y$ of index $r.$
\end{theorem1}

\begin{proof}
\normalfont Straightforward.
\end{proof}

\begin{definition1}
\normalfont Given two $g$-quasi metric spaces $(X,d_X)$ and $(Y,d_Y)$ of the same index $r,$ $d_{XY}$ is called the $g$-quasi metric product\index{product!g@$g$-quasi metric} of $d_X$ with $d_Y$ or simply product $g$-quasi metric on $X\times Y.$ 
\end{definition1}

Clearly if $(X,d_X)$ and $(Y,d_Y)$ are metric spaces, then $d_{XY}$ defines the product metric on $X\times Y$.

\begin{definition1}
\normalfont Let $(X,d_X)$ and $(Y,d_Y)$ be $g$-quasi metric spaces of indices $r_1$ and $r_2$ respectively. A mapping $f:X\to Y$ is said to be $g$-uniformly continuous\index{g@$g$-uniformly continuous mapping} if for $\epsilon>r_2$ there exists $\delta>r_1$ such that $d_X(x_1,x_2)<\delta\implies d_Y(f(x_1),f(x_2))<\epsilon,$ $\forall~x_1,x_2\in X.$
\end{definition1}

Clearly if $(X,d_X)$ and $(Y,d_Y)$ are metric spaces, then every $g$-uniformly continuous mapping from $X$ to $Y$ is uniformly continuous as a mapping between metric spaces.

It is known that if $(X,d)$ is a metric space, then the distance function $d:X\times X\to\mathbb R$ is uniformly continuous where $X\times X$ is equipped with the product metric and $\mathbb R$ with the usual metric. However given a $g$-quasi metric space $(X,d),$ the mapping $d:X\times X\to\mathbb R$ may not be $g$-uniformly continuous where $X\times X$ is equipped with the product metric and $\mathbb R$ with the usual metric (recall that, it is a $g$-quasi metric of index $0$). We consider the following example.

\begin{example1}
\normalfont Consider the $g$-quasi metric space $(\mathbb R,d),$ defined in Example \ref{gquasi_mainexm}, of index $r>0.$ We show that $d:(\mathbb R\times\mathbb R,d')\to(\mathbb R,d_u)$ is not $g$-uniformly continuous, where $d'$ is the $g$-quasi metric product of $d$ with itself on $\mathbb R\times\mathbb R$ and $d_u$ is the usual metric on $\mathbb R$.

Suppose otherwise. Then for $\epsilon=\frac{r}{2},$ $\exists~\delta>r$ such that $d'((x,y),(x',y'))<\delta\implies|d(x,y)-d(x',y')|<\epsilon,$ $\forall~(x,y),(x',y')\in\mathbb R\times\mathbb R.$

In particular, for $(x',y')=(0,0),$ $d'((x,y),(0,0))<\delta\implies|d(x,y)-d(0,0)|<\frac{r}{2},$ $\forall~(x,y)\in\mathbb R\times\mathbb R.$

That is, $\max\{d(x,0),d(y,0)\}<\delta\implies|d(x,y)-d(0,0)|<\frac{r}{2},$ $\forall~(x,y)\in\mathbb R\times\mathbb R.$

Choose $n\in\mathbb N\backslash\{1\}$ such that $\frac{\delta-r}{(n-1)^2}<\frac{r}{2}.$ Then $\frac{\delta-r}{n(n-1)}<\frac{r}{2}.$

Set $x=r+\frac{\delta-r}{n}$ and $y=r+\frac{\delta-r}{n-1}.$

Then $|x-0|=r+\frac{\delta-r}{n}>r$ and $|y-0|=r+\frac{\delta-r}{n-1}>r.$

Consequently, $d(x,0)=r+\frac{\delta-r}{n}<\delta$ and $d(y,0)=r+\frac{\delta-r}{n-1}<\delta$ whence, $\max\{d(x,0),$ $d(y,0)\}$ $<\delta.$

However, $|x-y|=\frac{\delta-r}{n(n-1)}<\frac{r}{2}\le r\implies d(x,y)=2r,$ a contradiction since $|d(x,y)-d(0,0)|<\frac{r}{2}.$

Hence $d:(\mathbb R\times\mathbb R,d')\to(\mathbb R,d_u)$ is not $g$-uniformly continuous.
\end{example1} 

\section{\textbf{Completeness, Lebesgue Property and (Weak) \texorpdfstring{$G$}--Completeness in \texorpdfstring{$g$}--Quasi Metric Spaces}}\label{section 5}

In this section, we extend the study of completeness, Lebesgue property and weak $G$-completeness for $g$-quasi metric spaces using the extended notion of Cauchy, $G$-Cauchy and pseudo-Cauchy sequences. 

\begin{definition1}
\normalfont Let $(x_n)$ be a sequence in a $g$-quasi metric space $(X,d)$ of index $r$ and $c\in X.$ Then

(i) $(x_n)$ is said to be convergent\index{convergent sequence} to $c$ in $(X,d)$ if it is so in $(X,\mu(d));$

(ii) $c$ is called a cluster point\index{cluster point} of $(x_n)$ in $(X,d)$ if it is so in $(X,\mu(d))$.

Clearly if $(x_n)$ is convergent to $c,$ then $c$ is a cluster point of $(x_n)$ (in $(X,d)$).
\end{definition1}

\begin{definition1}
\normalfont Let $(X,d)$ be a $g$-quasi metric space of index $r$ and $(x_n)$ be a sequence in $X.$ Then

(i) $(x_n)$ is called Cauchy\index{Cauchy sequence} if given $\epsilon>r$ there exists $k\in\mathbb N$ such that $d(x_m,x_n)<\epsilon,~\forall~m,n\ge k;$ 

(ii) $(x_n)$ is called $G$-Cauchy\index{g@$G$-Cauchy sequence!in $g$-quasi metric spaces} if given $\epsilon>r$ there exists $k\in\mathbb N$ such that $d(x_n,x_{n+1})<\epsilon,~\forall~n\ge k;$

(iii) $(x_n)$ is called pseudo-Cauchy\index{pseudo-Cauchy!sequence in $g$-quasi metric spaces} if given $\epsilon>r$ and $k\in\mathbb N$ there exist $p,q~(p\ne q)\in\mathbb N$ with $p,q\ge k$ such that $d(x_p,x_q)<\epsilon.$
\end{definition1}

\begin{definition1}
\normalfont A $g$-quasi metric space $(X,d)$ is said to be

(i) complete\index{complete!g@$g$-quasi metric space} if every Cauchy sequence converges to some point in it;

(ii) $G$-complete\index{g@$G$-complete!g@$g$-quasi metric space} if every $G$-Cauchy sequence converges to some point in it;

(iii) weak $G$-complete\index{weak $G$-complete!g@$g$-quasi metric space} if every $G$-Cauchy sequence has a cluster point in it;

(iv) Lebesgue\index{Lebesgue!g@$g$-quasi metric space} if every pseudo-Cauchy sequence having distinct terms has a cluster point in it;

(v) strongly Lebesgue\index{strongly Lebesgue!g@$g$-quasi metric space} if every pseudo-Cauchy sequence has a cluster point in it.
\end{definition1}

Clearly for $g$-quasi metric spaces  we have the following chain of implications:\\

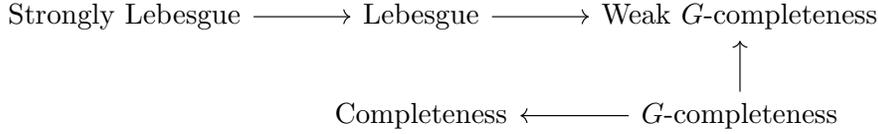
\noindent\begin{tikzcd}
\text{Strongly Lebesgue} \arrow{r} & \text{Lebesgue} \arrow{r} & \text{Weak }G\text{-completeness} \\
    & \text{Completeness} & \arrow{l} \arrow{u} G\text{-completeness}\\
\end{tikzcd}\\

In what follows, we show that for each $r>0,$ $(\mathbb R,d)$ of index $r,$ as defined in Example \ref{gquasi_mainexm}, is not weak $G$-complete.

\begin{example1}
\normalfont Consider the sequence $(x_n)$ in $(\mathbb R,d)$ where $x_n=rn-\frac{r}{n},~\forall~n\in\mathbb N.$

Then $\forall~n\in\mathbb N,$ $|x_{n+1}-x_n|=r+r\left(\frac{1}{n}-\frac{1}{n+1}\right)>r\implies d(x_{n+1},x_n)=r+r\left(\frac{1}{n}-\frac{1}{n+1}\right)\implies d(x_n,x_{n+1})=r+r\left(\frac{1}{n}-\frac{1}{n+1}\right).$ 

Choose $\epsilon>r.$ Then there is $k\in\mathbb N$ such that $\frac{r}{n(n+1)}<\epsilon-r,~\forall~n\ge k$ and hence, $d(x_n,x_{n+1})<\epsilon,~\forall~n\ge k.$ Thus $(x_n)$ is $G$-Cauchy in $(\mathbb R,d).$

If possible, let there is a cluster point $c$ of $(x_n)$ in $(\mathbb R,d).$ 

Since $(x_n)$ is a sequence of distinct terms, $B(c,\frac{3r}{2})$ contains infinitely many elements of $(x_n).$

However $B(c,\frac{3r}{2})=\{y\in\mathbb R:d(c,y)<\frac{3r}{2}\}=\left(c-\frac{3r}{2},c-r\right)\bigcup\left(c+r,c+\frac{3r}{2}\right)\bigcup\{c\}$ which clearly contains finitely many elements of $(x_n),$ a contradiction.

Hence $(\mathbb R,d)$ is not weak $G$-complete.
\end{example1}

\begin{lemma1}\label{ch6_lem_cauchy}
\normalfont Let $(X,d_X)$ and $(Y,d_Y)$ be $g$-quasi metric spaces of the same index $r.$ A sequence $((x_n,y_n))$ is Cauchy in $(X\times Y,d_{XY})$ if and only if $(x_n)$ and $(y_n)$ are Cauchy in $(X,d_X)$ and $(Y,d_Y)$ respectively.
\end{lemma1}

\begin{proof}
\normalfont Let $((x_n,y_n))$ be Cauchy in $(X\times Y,d_{XY}).$ Choose $\epsilon>r.$ Then $\exists~k\in\mathbb N$ such that $d_{XY}((x_m,y_m),(x_n,y_n))<\epsilon,~\forall~m,n\ge k.$ That is, $d_{X}(x_m,x_n),d_{Y}(y_m,y_n)<\epsilon,~\forall~m,n\ge k.$ Then $(x_n)$ and $(y_n)$ are Cauchy in $(X,d_X)$ and $(Y,d_Y)$ respectively.

\textit{Conversely}, let $(x_n)$ and $(y_n)$ be Cauchy in $(X,d_X)$ and $(Y,d_Y)$ respectively. Choose $\epsilon>r.$ Then $\exists~p,q\in\mathbb N$ such that $d_X(x_m,x_n)<\epsilon,~\forall~m,n\ge p$ and $d_Y(y_m,y_n)<\epsilon,~\forall~m,n\ge q.$

Set $r=\max\{p,q\}.$ Then $d_{XY}((x_m,y_m),(x_n,y_n))<\epsilon,~\forall~m,n\ge r.$

Hence $((x_n,y_n))$ is Cauchy in $(X\times Y,d_{XY}).$
\end{proof}

Similar chain of arguments yield the following results that we state without proof.

\begin{lemma1}
\normalfont Let $(X,d_X)$ and $(Y,d_Y)$ be $g$-quasi metric spaces of the same index $r.$ A sequence $((x_n,y_n))$ is $G$-Cauchy in $(X\times Y,d_{XY})$ if and only if $(x_n)$ and $(y_n)$ are $G$-Cauchy in $(X,d_X)$ and $(Y,d_Y)$ respectively.
\end{lemma1}

\begin{lemma1}\label{ch6_(xn,yn)pseudo-Cauchy}
\normalfont Let $(X,d_X)$ and $(Y,d_Y)$ be $g$-quasi metric spaces of the same index $r.$ If $((x_n,y_n))$ is a pseudo-Cauchy sequence in $(X\times Y,d_{XY})$ then $(x_n)$ and $(y_n)$ are pseudo-Cauchy in $(X,d_X)$ and $(Y,d_Y)$ respectively.
\end{lemma1}

Then converse of Lemma \ref{ch6_(xn,yn)pseudo-Cauchy} is not true. In support, we produce the following example.

\begin{example1}
\normalfont Consider the $g$-quasi metric space $(\mathbb R,d),$ as defined in Example \ref{gquasi_mainexm}, of index $1.$ Then the sequences $(x_n)$ and $(y_n)$ in $\mathbb R,$ defined by $$x_n=\begin{cases}1 &\text{if } n\text{ is odd}\\10^n &\text{if } n\text{ is even}\end{cases}$$ and $$y_n=\begin{cases}10^n &\text{if } n\text{ is odd}\\1 &\text{if } n\text{ is even}\end{cases}$$ $\forall~n\in\mathbb N,$ are pseudo-Cauchy in $(\mathbb R,d).$ 

However $((x_n,y_n))$ is not pseudo-Cauchy in $\mathbb R\times\mathbb R$ (where $\mathbb R\times\mathbb R$ is equipped with $d',$ the $g$-quasi metric product of $d$ with itself). In fact, for any pair of positive integers $m,q~(m\ne q)$ with $m,q\ge 1,$ we obtain $d'((x_m,y_m),(x_q,y_q))>2,$ by considering even and odd cases separately for $m$ and $q.$ Thus $((x_n,y_n))$ is not pseudo-Cauchy in $(\mathbb R\times\mathbb R,d').$
\end{example1}

\begin{theorem1}
\normalfont Let $(X,d_X)$ and $(Y,d_Y)$ be $g$-quasi metric spaces of the same index $r.$ Then $(X\times Y,d_{XY})$ is complete if and only if $(X,d_X)$ and $(Y,d_Y)$ are complete. 
\end{theorem1}

\begin{proof}
\normalfont Let $(X\times Y,d_{XY})$ be complete.

We first show that $(X,d_X)$ is complete. Choose a Cauchy sequence $(x_n)$ in $(X,d_X)$ and fix $y\in Y.$ Then by Lemma \ref{ch6_lem_cauchy}, $((x_n,y))$ is Cauchy in $(X\times Y,d_{XY}).$ Since $(X\times Y,d_{XY})$ is complete, $((x_n,y))$ converges to a point $(a,b)$ in $(X\times Y,d_{XY}).$

We claim that $(x_n)$ is convergent to $a$ in $(X,d_{X}).$

Let $V$ be a generalized open set in $(X,\mu(d_X))$ containing $a.$ Then there exist $p\in X,\delta>r$ such that $a\in B_{d_X}(p,\delta)\subset V.$

Since $B_{d_{XY}}((p,b),\delta)$ is an open set in $(X\times Y,\mu(d_{XY}))$ containing $(a,b),~\exists~k\in\mathbb N$ such that $(x_n,y)\in B_{d_{XY}}((p,b),\delta),$ $\forall~n\ge k.$ 

Thus $\max\{d_X(p,x_n),d_Y(b,y)\}<\delta,~\forall~n\ge k\implies d_X(p,x_n)<\delta,~\forall~n\ge k.$

i.e., $x_n\in B_{d_X}(p,\delta)\subset V,~\forall~n\ge k.$ Hence $(x_n)$ converges to $a$ in $(X,d_X)$ and so, $(X,d_X)$ is complete. Similarly $(Y,d_Y)$ is complete.

\textit{Conversely}, let $(X,d_X)$ and $(Y,d_Y)$ be complete.

Choose a Cauchy sequence $((x_n,y_n))$ in $(X\times Y,d_{XY}).$ By Lemma \ref{ch6_lem_cauchy}, $(x_n)$ and $(y_n)$ are Cauchy in $(X,d_X)$ and $(Y,d_Y)$ respectively. So by hypothesis, there exist $a\in X,b\in Y$ such that $(x_n)$ converges to $a$ in $(X,d_X)$ and $(y_n)$ to $b$ in $(Y,d_Y).$

We show that $((x_n,y_n))$ converges to $(a,b)$ in $(X\times Y,d_{XY}).$

Let $W$ be a generalized open set in $(X\times Y,\mu(d_{XY}))$ containing $(a,b).$ Then there exist $(p,q)\in X\times Y,\delta>r$ such that $(a,b)\in B_{d_{XY}}((p,q),\delta)\subset W.$ Consequently $a\in B_{d_X}(p,\delta)$ and $b\in B_{d_Y}(q,\delta).$

Since $(x_n)$ converges to $a$ and $(y_n)$ to $b,$ there exist $k_1,k_2\in\mathbb N$ such that $x_n\in B_{d_X}(p,\delta),~\forall~n\ge k_1$ and $y_n\in B_{d_Y}(q,\delta),~\forall~n\ge k_2.$

Set $k=\max\{k_1,k_2\}.$ Then $(x_n,y_n)\in B_{d_{XY}}((p,q),\delta)\subset W,~\forall~n\ge k.$ Thus $((x_n,y_n))$ converges to $(a,b)$ in $(X\times Y,d_{XY}).$ 

Thus $(X\times Y,d_{XY})$ is complete.
\end{proof}

Similar chain of arguments yield the following results that we state without
proof.

\begin{theorem1}
\normalfont Let $(X,d_X)$ and $(Y,d_Y)$ be $g$-quasi metric spaces of the same index $r.$ Then $(X\times Y,d_{XY})$ is $G$-complete if and only if $(X,d_X)$ and $(Y,d_Y)$ are $G$-complete. 
\end{theorem1}

\begin{theorem1}
\normalfont Let $(X,d_X)$ and $(Y,d_Y)$ be $g$-quasi metric spaces of the same index $r.$ If $(X\times Y,d_{XY})$ is weak $G$-complete then $(X,d_X)$ and $(Y,d_Y)$ are weak $G$-complete. 
\end{theorem1}

\begin{theorem1}
\normalfont Let $(X,d_X)$ and $(Y,d_Y)$ be $g$-quasi metric spaces of the same index $r.$ If $(X\times Y,d_{XY})$ is (strongly) Lebesgue then $(X,d_X)$ and $(Y,d_Y)$ are (strongly) Lebesgue. 
\end{theorem1}


\begin{thebibliography}{99}

\bibitem{1sa-g-com} Adhya, S., and Deb Ray, A. (2022). On weak $G$-completeness for fuzzy metric spaces. Soft Computing, 26(5), 2099-2105.

\bibitem{l1} Atsuji, M. (1958). Uniform continuity of continuous functions of metric spaces. Pacific J. Math., 8(4), 11--16.

\bibitem{csaszar2004extremally} Cs{\'a}sz{\'a}r, {\'A}. (2004). Extremally disconnected generalized topologies. Annales Univ. Sci. Budapest, 47, 91--96.

\bibitem{csaszar1997generalized} Cs{\'a}sz{\'a}r, {\'A}. (1997). Generalized open sets. Acta Math. Hungar. 75(1--2), 65--87. 

\bibitem{czaszar1} Cs{\'a}sz{\'a}r, {\'A}. (2002). Generalized topology, generalized continuity. Acta Math. Hungar., 96(4), 351--357.

\bibitem{csaszar2007modification} Cs{\'a}sz{\'a}r, {\'A}. (2007). Modification of generalized topologies via hereditary classes. Acta Math. Hungar., 115(1--2), 29--36.

\bibitem{czaszarproduct} Cs{\'a}sz{\'a}r, {\'A}. (2009). Product of generalized topologies. Acta Math. Hungar., 123(1--2), 127--132.

\bibitem{raybhowmik2015} Deb Ray, A. and Bhowmick, R. (2015). On generalized quasi-uniformity and generalized quasi-uniformizable supratopological spaces. J. Adv. Stud. Top., 6(2), 74--81.

\bibitem{2greg-g-com} Gregori, V., Miñana, J. J., Roig, B., and Sapena, A. (2018). On completeness in metric spaces and fixed point theorems. Results in Mathematics, 73, 1--13.

\bibitem{3greg2g-com} Gregori, V., Miñana, J. J., and Sapena, A. (2018). On Banach contraction principles in fuzzy metric spaces. Fixed point theory, 19(1), 235--247.

\bibitem{l2} Kundu, S., and Jain, T. (2006). Atsuji spaces: equivalent conditions. Topology Proc., 30(1), 301--325.

\bibitem{levine1963semi} Levine, N. (1963). Semi-open sets and semi-continuity in topological spaces. Amer. Math. Monthly, 70, 36--41.

\bibitem{mashhour1983supratopological} Mahmoud, F. S. and Khedr, F. H. and Mashhour, A. S. and Allam, A. A. (1983). On supratopological spaces. Indian J. Pure Appl. Math., 14(4), 502--510.

\bibitem{abd1982precontinuous} Mashhour, A. S., Abd El-Monsef, M. E. and El-Deeb, S. N. (1982). On precontinuous and weak precontinuous mappings. Proc. Math. Phys. Soc. Egypt, 53, 47--53.

\bibitem{Njstad1965OnSC} Olav N. (1965). On some classes of nearly open sets. Pacific J. Math., 15(3), 961--970.

\bibitem{l3} Toader, G. (1978). On a problem of Nagata. Mathematica (Cluj), 20(43), 78-79.

\bibitem{wilson1931quasi} Wilson, W. A. (1931). On quasi-metric spaces. Amer. J. Math., 53(3), 675--684.

\bibitem{wu-zhu} Wu, X. and Zhu, P. (2013). Generalized product topology. Commun. Korean Math. Soc., 28(4), 819--825.

\bibitem{sarsak_weak_sep} Xun, GE and Ying, GE (2010). $\mu$-{S}eparations in generalized topological spaces. Appl. Math. J. Chinese Univ., 25(2), 243--252.

\end{thebibliography}
\end{document}